\documentclass[letterpaper,11pt]{article}
\usepackage[margin=1in]{geometry}

\usepackage{amsmath, amsthm, amssymb}
\usepackage{enumerate}
\usepackage{cite}
\usepackage{mathtools}
\usepackage{cancel}
\usepackage{siunitx}
\usepackage{tikz}
\usetikzlibrary{matrix,calc}
\usepackage{pgfplots}
\usepackage{url}
\pgfplotsset{compat=1.12}

\newtheorem{theo}{Theorem}[section]

\newtheorem{lemma}[theo]{Lemma}

\tikzset{square matrix/.style={
    matrix of nodes,
    column sep=-\pgflinewidth, row sep=-\pgflinewidth,
    nodes={draw,
      minimum height=15pt,
      anchor=center,
      text width=13pt,
      align=center,
      inner sep=0pt
    },
  },
  square matrix/.default=2cm
}

\newcommand\BibTeX{{\rmfamily B\kern-.05em \textsc{i\kern-.025em b}\kern-.08em
T\kern-.1667em\lower.7ex\hbox{E}\kern-.125emX}}

\newcommand{\cp}{\mathcal P}

\begin{document}

\title{Computing autotopism groups of partial Latin rectangles: a pilot study}

\author{
Rebecca J.\ Stones\\ \small College of Computer Science, Nankai University, Tianjin, China.\\
\small \url{rebecca.stones82@gmail.com}\\ \\
Ra\'{u}l M.\ Falc\'{o}n\\
\small  School of Building Engineering, University of Seville, Spain.\\ \small \url{rafalgan@us.es}\\ \\
Daniel Kotlar\\
\small  Department of Computer Science, Tel-Hai College, Upper Galilee, Israel.\\ \small \url{dannykotlar@gmail.com}\\ \\
Trent G. Marbach\\
\small  College of Computer Science, Nankai University, Tianjin, China.\\ \small \url{trent.marbach@outlook.com}\\ \\
}
\date{\today}

\maketitle

\begin{abstract}
Computing the autotopism group of a partial Latin rectangle can be performed in a variety of ways.  This pilot study has two aims: (a) to compare these methods experimentally,  and (b) to identify the design goals one should have in mind for developing practical software. To this end, we compare six families of algorithms (two backtracking methods and four graph automorphism methods), with and without the use of entry invariants, on two test suites.  We consider two entry invariants: one determined by the frequencies of row, column, and symbol representatives, and one determined by $2 \times 2$ submatrices. We find: (a) with very few entries, many symmetries often exist, and these should be identified mathematically rather than computationally, (b) with an intermediate number of entries, a quick-to-compute entry invariant was effective at reducing the need for computation, (c) with an almost-full partial Latin rectangle, more sophisticated entry invariants are needed, and (d) the performance for (full) Latin squares is significantly poorer than other partial Latin rectangles of comparable size, obstructed by the existence of Latin squares with large (possibly transitive) autotopism groups.
\end{abstract}


\section{Introduction}

Let $[n]=\{1,2,\ldots,n\}$.  An $r \times s$ \emph{partial Latin rectangle} is an $r \times s$ array containing symbols from $[n] \cup \{\cdot\}$, where  each symbol in $[n]$ appears at most once within each row and each column. If $r=s=n$, then this is a \emph{partial Latin square of order $n$}.

Let $\mathrm{PLR}(r,s,n)$ denote the set of $r \times s$ partial Latin rectangles. Every $L\in\mathrm{PLR}(r,s,n)$ is uniquely determined by its \emph{entry set} $$\mathrm{Ent}(L):=\{(i,j,L[i,j]):i\in [r], j\in [s], L[i,j] \in [n]\}.$$
Any triple $(i,j,L[i,j]) \in \mathrm{Ent}(L)$ is an \emph{entry} of $L$, while any pair $(i,j) \in [r] \times [s]$ is a \emph{cell}. The latter is considered \emph{empty} if it contains the symbol $\cdot$, in which case, $L[i,j]$ is \emph{undefined}.  If $L \in \mathrm{PLR}(r,s,n)$ does not have empty cells, then it is a \emph{Latin rectangle} (or a \emph{Latin square}, if $r=s=n$).  This constitutes a generalization of the standard definition of a Latin rectangle \cite{Stones2009b}, which requires that $s=n$.

Let $S_m$ denote the symmetric group on $[m]$. A triple of permutations $\theta:=(\alpha,\beta,\gamma) \in S_r \times S_s \times S_n$ acts on $\mathrm{PLR}(r,s,n)$, with $\theta$ acting on $L \in \mathrm{PLR}(r,s,n)$ by permuting its rows, columns, and symbols by $\alpha$, $\beta$, and $\gamma$, respectively. The triple $\theta$ is called an \emph{isotopism} and both $L$ and the resulting partial Latin rectangle $L^\theta\in \mathrm{PLR}(r,s,n)$ are \emph{isotopic}. Moreover, 
$$\mathrm{Ent}(L^{\theta})=\left\{\big(\alpha(i),\beta(j),\gamma(L[i,j])\big)\colon\, \big(i,j,L[i,j]\big)\in \mathrm{Ent}(L)\right\},$$
where, from here on, we use $(i,j,L[i,j])^\theta$ to denote $(\alpha(i),\beta(j),\gamma(L[i,j]))$.  If $L=L^\theta$, then $\theta$ is said to be an \emph{autotopism} of $L$.  The set $\mathrm{Atop}(L)$ of autotopisms of $L$ forms a group, called the \emph{autotopism group} of $L$.

Autotopisms of (partial) Latin squares have arisen in a range of topics and are also studied in their own right \cite{Artzy1954, Bailey1982, Browning2013, BryantBuchananWanless2009, CavenaghStones2010a, Drisko1997b, Falcon2009, FalconMartinMorales2007, FalconNunez2007, IhrigIhrig2008, KerbySmith2009, KerbySmith2010, kot12, McKayWanlessZhang2015, Mendis2016, MengZhengZheng2008, StonesVojtechovskyWanless2012, WanlessIhrig2005}. In particular, they are beginning to find applications in cryptography \cite{Ganfornina2006, KongEtAl, StonesSuLiuWangLin}, erasure codes \cite{YiEtAl,StonesStrongKPlex} (relating to Blackburn partial Latin squares \cite{Wanless2004b}), and graph colouring games \cite{FalconAndres2018a, FalconAndres2018b}.  However, autotopisms of (partial) Latin rectangles, have only briefly been studied in recent papers \cite{Falcon2013, FalconStones2015, FalSto16}. 

The standard way of computing autotopism groups of Latin squares was given by McKay, Meynert, and Myrvold \cite{McKayMeynertMyrvold2007}, who describe a graph whose automorphism group is isomorphic to the autotopism group of the corresponding Latin square.  A generalization of this method has already been used for studying partial Latin squares with trivial symmetry groups \cite{Stones13}. For autotopisms of Latin squares, Kotlar \cite{kot14} gave a method which enumerates cycles belonging to entries, which can be thought of as a special case of, what we call in this paper, row invariants.

Almost all Latin squares have trivial autotopism groups \cite{Cameron2015, McKayWanless2005}, so, in almost all cases, computing such a group boils down to eliminating the possible existence of a non-trivial autotopism.  This procedure cannot be, however,  straightforwardly generalized to (partial) Latin rectangles, because: (a) two empty rows (or two empty columns) in a partial Latin rectangle can be swapped, giving rise to a large family of partial Latin rectangles with a non-trivial autotopism; (b) $2 \times n$ Latin rectangles always have non-trivial symmetries; and (c) partial Latin squares with very few entries must have non-trivial symmetries \cite{Stones13}.

In this paper, six families of algorithms for computing this group are implemented and experimentally compared.  Some of them have been used in prior research on partial Latin rectangles, while others are introduced here.

\section{Autotopisms of partial Latin rectangles}

\subsection{Introduction}

We start with a characterization of the existence of an autotopism of a partial Latin rectangle with specified row and column permutations, but with an unspecified symbol permutation. We defer proofs to a subsequent paper.

\begin{lemma}\label{lm:atop_cond}
Let $L \in \mathrm{PLR}(r,s,n)$.  If $\alpha \in S_r$ and $\beta \in S_s$, then there exists $\gamma \in S_n$ for which $\theta=(\alpha,\beta,\gamma)\in\mathrm{Atop}(L)$ if and only if
\begin{enumerate}
 \item\label{cond:defined} $L[\alpha(i),\beta(j)]$ is defined for all $(i,j,L[i,j]) \in \mathrm{Ent}(L)$, and
 \item\label{cond:symbols} $L[\alpha(i),\beta(j)]=L[\alpha(i'),\beta(j')]$, whenever $L[i,j]=L[i',j'] \in [n]$.
\end{enumerate}
The permutations $\gamma$ for which $\theta\in\mathrm{Atop}(L)$ are those satisfying $\gamma(L[i,j])=L[\alpha(i),\beta(j)]$, for all $(i,j,L[i,j]) \in \mathrm{Ent}(L)$. \hfill \qed
\end{lemma}

Lemma~\ref{lm:atop_cond} suggests an algorithm for simultaneously determining whether or not a row permutation $\alpha$ and a column permutation $\beta$ participate in an autotopism of a partial Latin rectangle $L \in \mathrm{PLR}(r,s,n)$, and identifying which symbol permutations $\gamma$ are possible.  We proceed entry by entry, setting $\gamma(L[i,j])=L[\alpha(i),\beta(j)]$ for every entry $(i,j,L[i,j])\in\mathrm{Ent}(L)$.  Three clashes can arise: (a) $L[\alpha(i),\beta(j)]$ is undefined; (b) $\gamma(L[i,j])$ was previously set to something else; or (c) $\gamma^{-1}(L[\alpha(i),\beta(j)])$ was previously set to something else.  If no clashes arise, then the triple $(\alpha,\beta,\gamma)$ is an autotopism for all completions of $\gamma$ to a permutation in $S_n$.

Lemma~\ref{lm:atop_cond} also enables one to compute the autotopism group while ignoring empty rows, empty columns, and unused symbols. Specifically, suppose $L$ has $r'\leq r$ non-empty rows, $s'\leq s$ non-empty columns and $n'\leq n$ used symbols. Let $\overline{L}\in \mathrm{PLR}(r',s',n')$ denote the partial Latin rectangle that results after eliminating the empty rows and columns of $L$ and relabeling the used symbols of $L$ to $[n']$. Thus, the following result holds.

\begin{lemma}\label{lm:atop_card}
$|\mathrm{Atop}(L)|=(r-r')!\, (s-s')!\, (n-n')!\, |\mathrm{Atop}(\overline{L})|$. \hfill \qed
\end{lemma}

\subsection{Entry invariants}

An \emph{entry invariant} $\mathcal{P}$ is any property of the entries of partial Latin rectangles which is preserved under autotopisms. That is, if $L$ is a partial Latin rectangle, then $\mathcal{P}(e)=\mathcal{P}(e^\theta)$, for all $e\in \mathrm{Ent}(L)$ and $\theta\in\mathrm{Atop}(L)$. We also define a \emph{row invariant}, \emph{column invariant}, and \emph{symbol invariant} of $L$ to be respectively any property of its rows, columns, or symbols which is preserved under autotopisms. Given an entry invariant $\mathcal{P}$, a useful row invariant is the multiset of entry invariants of the entries in that row, denoted $\mathcal{P}_R$, and likewise for columns, denoted $\mathcal{P}_C$. The symbol invariant $\cp_S$ is the multiset of entry invariants of the entries containing the particular symbol in the third component.

\begin{lemma}\label{lm:atop_inv2}  Let $L \in \mathrm{PLR}(r,s,n)$ have entry invariant $\mathcal{P}$. Then, $L$ has a trivial autotopism group if one of the following conditions holds.
\begin{enumerate}[a)]
\item $\mathcal{P}(e)\neq \mathcal{P}(e')$, for every pair of distinct entries $e,e'\in\mathrm{Ent}(L)$.
\item There exist distinct $X,Y\in\{R,C,S\}$ such that $\mathcal{P}_X(i)\neq \mathcal{P}_X(i')$ whenever $i\ne i'$, and $\mathcal{P}_Y(j)\neq \mathcal{P}_Y(j')$ whenever $j\ne j'$, where $i,i',j$ and $j'$ are taken from the appropriate sets of rows, columns, or symbols.  \hfill \qed
\end{enumerate}
\end{lemma}

Every entry invariant $\mathcal{P}$ of $L$ determines a partition of its set $\mathrm{Ent}(L)$ whereby we replace all the entries on $L$ by the index of their entry invariant (we subsequently give examples in \eqref{eq:strongentryinvariants} and \eqref{eq:squareinvariants}). Let $\mathbf{c}^{\mathcal{P}}_j(L)$  denote the $j$-th column of the resulting new array, for all $j\in [s]$. Let $\overline{\mathbf{c}}_j(L) \subseteq [r]$ denote the set of row indices $i$ for which the cell $(i,j)$ is non-empty.

\begin{lemma}\label{lm:col_vec}
Let $L \in \mathrm{PLR}(r,s,n)$ and $\theta=(\alpha,\beta,\gamma)\in\mathrm{Atop}(L)$. If $\mathcal{P}$ is invariant under symbol permutations, then
\begin{enumerate}[a)]
\item $\mathbf{c}^{\mathcal{P}}_j(L^{(\alpha,\mathrm{id},\mathrm{id})})[i]=\mathbf{c}^{\mathcal{P}}_{\beta(j)}(L)[i]$, for all $i\in [r]$ and $j\in [s]$.
\item $\overline{\mathbf{c}}_j(L)=\overline{\mathbf{c}}_{\beta(j)}(L)$, for all $j\in [s]$.
% \item $\delta^{\mathcal{P}}_j(L) = \delta^{\mathcal{P}}_{\beta(j)}(L)$, for all $j\in [s]$.  \hfill \qed
\end{enumerate}
\end{lemma}

In this paper, all entry invariants are invariant under symbol permutations.

\subsubsection{Strong entry invariants}\label{SEI}

The \emph{strong entry invariant} \cite{FalSto16} of an entry $(i,j,L[i,j])$ is the triple $(a,b,c)$ formed by the number $a$ of entries in row $i$, the number $b$ of entries in column $j$, and the number $c$ of times the symbol $L[i,j]$ appears in $L$.  For example:
\begin{equation}\label{eq:strongentryinvariants}
\begin{array}{c}
\begin{array}{c}
\begin{tikzpicture}
\matrix[square matrix]{
|[fill=black!10]| 1 & |[fill=white]| $\cdot$ & |[fill=black!21]| 2 & |[fill=white]| $\cdot$ & 
|[fill=white]| $\cdot$ & |[fill=white]| $\cdot$ & |[fill=black!33]| 3 & |[fill=white]| $\cdot$ & 
|[fill=white]| $\cdot$ \\
|[fill=black!21]| 2 & |[fill=white]| $\cdot$ & |[fill=white]| $\cdot$ & |[fill=black!45]| 4 & 
|[fill=black!10]| 1 & |[fill=black!56]| 5 & |[fill=black!68]| 6 & |[fill=white]| $\cdot$ & 
|[fill=black!80]| 7 \\
|[fill=white]| $\cdot$ & |[fill=black!10]| 1 & |[fill=black!56]| 5 & |[fill=black!33]| 3 & 
|[fill=white]| $\cdot$ & |[fill=black!45]| 4 & |[fill=white]| $\cdot$ & |[fill=white]| $\cdot$ & 
|[fill=white]| $\cdot$ \\
|[fill=white]| $\cdot$ & |[fill=black!21]| 2 & |[fill=white]| $\cdot$ & |[fill=black!56]| 5 & 
|[fill=white]| $\cdot$ & |[fill=black!33]| 3 & |[fill=white]| $\cdot$ & |[fill=black!45]| 4 & 
|[fill=white]| $\cdot$ \\
|[fill=black!45]| 4 & |[fill=black!33]| 3 & |[fill=white]| $\cdot$ & |[fill=white]| $\cdot$ & 
|[fill=black!56]| 5 & |[fill=white]| $\cdot$ & |[fill=black!10]| 1 & |[fill=white]| $\cdot$ & 
|[fill=black!21]| 2 \\
|[fill=white]| $\cdot$ & |[fill=white]| $\cdot$ & |[fill=white]| $\cdot$ & |[fill=white]| $\cdot$ & 
|[fill=black!21]| 2 & |[fill=white]| $\cdot$ & |[fill=white]| $\cdot$ & |[fill=black!10]| 1 & 
|[fill=black!33]| 3 \\
};
\end{tikzpicture}
\end{array}
\xrightarrow{\substack{\text{strong entry invariants} \\ \text{relabeled $1,2,\ldots$}}}
\begin{array}{c}
\begin{tikzpicture}
\matrix[square matrix]{
|[fill=black!10]| 1 & |[fill=white]| $\cdot$ & |[fill=black!17]| 2 & |[fill=white]| $\cdot$ & 
|[fill=white]| $\cdot$ & |[fill=white]| $\cdot$ & |[fill=black!10]| 1 & |[fill=white]| $\cdot$ & 
|[fill=white]| $\cdot$ \\
|[fill=black!25]| 3 & |[fill=white]| $\cdot$ & |[fill=white]| $\cdot$ & |[fill=black!33]| 4 & 
|[fill=black!25]| 3 & |[fill=black!33]| 4 & |[fill=black!41]| 5 & |[fill=white]| $\cdot$ & 
|[fill=black!41]| 5 \\
|[fill=white]| $\cdot$ & |[fill=black!48]| 6 & |[fill=black!56]| 7 & |[fill=black!48]| 6 & 
|[fill=white]| $\cdot$ & |[fill=black!64]| 8 & |[fill=white]| $\cdot$ & |[fill=white]| $\cdot$ & 
|[fill=white]| $\cdot$ \\
|[fill=white]| $\cdot$ & |[fill=black!48]| 6 & |[fill=white]| $\cdot$ & |[fill=black!64]| 8 & 
|[fill=white]| $\cdot$ & |[fill=black!48]| 6 & |[fill=white]| $\cdot$ & |[fill=black!56]| 7 & 
|[fill=white]| $\cdot$ \\
|[fill=black!72]| 9 & |[fill=black!80]| 10 & |[fill=white]| $\cdot$ & |[fill=white]| $\cdot$ & 
|[fill=black!72]| 9 & |[fill=white]| $\cdot$ & |[fill=black!80]| 10 & |[fill=white]| $\cdot$ & 
|[fill=black!80]| 10 \\
|[fill=white]| $\cdot$ & |[fill=white]| $\cdot$ & |[fill=white]| $\cdot$ & |[fill=white]| $\cdot$ & 
|[fill=black!10]| 1 & |[fill=white]| $\cdot$ & |[fill=white]| $\cdot$ & |[fill=black!17]| 2 & 
|[fill=black!10]| 1 \\
};
\end{tikzpicture}
\end{array}
\end{array}
\end{equation}
The multisets of strong entry invariants in each row and column imply that in any autotopism $(\alpha,\beta,\gamma)$ of this partial Latin rectangle, $\alpha$ and $\beta$ fix the partitions $\{ \{ 1, 6 \}, \{ 2 \}, \{ 3, 4 \}, \{ 5 \} \}$ and $\{ \{ 1, 5 \}, \{ 2 \}, \{ 3, 8 \}, \{ 4, 6 \}, \{ 7, 9 \} \}$ set-wise, respectively.  It turns out there is a single non-trivial autotopism with row permutation $(16)(34)$, column permutation $(15)(38)(46)(79)$ and symbol permutation $(12)(45)(67)$. Thus, no finer partitions of the sets $[6]$ and $[9]$ can be achieved using row and column invariants, respectively.

\subsubsection{Square invariants}\label{square}

Given a partial Latin rectangle $L\in \mathrm{PLR}(r,s,n)$, any entry $(i,j,k)$ belongs to exactly $(r-1)(s-1)$ distinct $2 \times 2$ submatrices of $L$, a general one having the following form:
$$
\begin{array}{c|cc}
   & j & j' \\
\hline
i  & k & x \\
i' & y & z \\
\end{array}
$$
which may have some of the following five properties: (a) $x$ is undefined, (b) $y$ is undefined, (c) $z$ is undefined, (d) $k=z$, and (e) $x=y$.  This gives a maximum of $2^5=32$ possibilities, whose enumeration gives a length-$32$ vector that sums to $(r-1)(s-1)$. This vector constitutes an entry invariant, which we call the \emph{square invariant}. Unlike strong entry invariants, the square invariants for Latin squares are not necessarily equal.  For example:
\begin{equation}\label{eq:squareinvariants}
\begin{array}{c}
\begin{tikzpicture}
\matrix[square matrix]{
|[fill=black!27]| 2 & |[fill=black!10]| 1 & |[fill=black!45]| 3 & |[fill=black!62]| 4 & |[fill=black!80]| 5 \\
|[fill=black!10]| 1 & |[fill=black!62]| 4 & |[fill=black!27]| 2 & |[fill=black!80]| 5 & |[fill=black!45]| 3 \\
|[fill=black!62]| 4 & |[fill=black!45]| 3 & |[fill=black!80]| 5 & |[fill=black!10]| 1 & |[fill=black!27]| 2 \\
|[fill=black!80]| 5 & |[fill=black!27]| 2 & |[fill=black!10]| 1 & |[fill=black!45]| 3 & |[fill=black!62]| 4 \\
|[fill=black!45]| 3 & |[fill=black!80]| 5 & |[fill=black!62]| 4 & |[fill=black!27]| 2 & |[fill=black!10]| 1 \\
};
\end{tikzpicture}
\end{array}
\xrightarrow{\substack{\text{$2 \times 2$ entry invariants} \\ \text{relabeled $1,2,\ldots$}}}
\begin{array}{c}
\begin{tikzpicture}
\matrix[square matrix]{
|[fill=black!10]| 1 & |[fill=black!45]| 2 & |[fill=black!10]| 1 & |[fill=black!10]| 1 & |[fill=black!45]| 2 \\
|[fill=black!45]| 2 & |[fill=black!10]| 1 & |[fill=black!10]| 1 & |[fill=black!10]| 1 & |[fill=black!45]| 2 \\
|[fill=black!10]| 1 & |[fill=black!10]| 1 & |[fill=black!10]| 1 & |[fill=black!45]| 2 & |[fill=black!45]| 2 \\
|[fill=black!10]| 1 & |[fill=black!10]| 1 & |[fill=black!45]| 2 & |[fill=black!10]| 1 & |[fill=black!45]| 2 \\
|[fill=black!45]| 2 & |[fill=black!45]| 2 & |[fill=black!45]| 2 & |[fill=black!45]| 2 & |[fill=black!80]| 3 \\
};
\end{tikzpicture}
\end{array}
\end{equation}

\section{Computing the autotopism group}

Here, we introduce the methods we consider for computing the autotopism group.

\subsection{Backtracking methods}\label{se:backtrack}

\subsubsection{Alpha-beta backtracking}

We start with the outputs of $\alpha$, $\beta$, and $\gamma$ as undefined, and use backtracking to determine the possible row permutations $\alpha$. At each level of the search tree we designate $$i \xmapsto{\alpha} a$$ for some $i,a \in [r]$ provided $\alpha^{-1}(a)$ is not already defined as something other than $i$.  Once $\alpha$ has been decided, we again use backtracking to establish the possible column permutations $\beta$.  At each level, we designate $$j \xmapsto{\beta} b$$ for some $j,b \in [s]$ provided $\beta^{-1}(b)$ is not already defined as something other than $j$.  Finally, for each pair $(\alpha,\beta)$, Lemma~\ref{lm:atop_cond} describes how to determine whether an autotopism $(\alpha,\beta,\gamma)$ exists, and how to find the possible symbol permutations $\gamma$.

In practice, both alpha- and beta-search trees can benefit from invariants, e.g., if we decide $\alpha(i)=a$, then rows $i$ and $a$ must have the same value for the invariants described in \ref{SEI} and \ref{square}.  In this regard, we can use Lemma~\ref{lm:col_vec} to prune the backtracking method described above. We refer to this improvement of the backtracking method as the \emph{CV method} (Column Vector). 

\subsubsection{Entrywise backtracking}

We start with the outputs of $\alpha$, $\beta$, and $\gamma$ as undefined, and at each level of the search tree, we designate $$(i,j,L[i,j]) \xmapsto{\theta} (a,b,L[a,b])$$ or equivalently, we simultaneously designate
\begin{align*}
i \xmapsto{\alpha} a, \qquad j \xmapsto{\beta} b, \qquad \text{and} \qquad L[i,j] \xmapsto{\gamma} L[a,b],
\end{align*}
provided none of the next six clashes arises: (a) $\alpha(i)$ is already defined as something other than $a$; (b) $\beta(j)$ is already defined as something other than $b$; (c) $\gamma(L[i,j])$ is already defined as something other than $L[a,b]$; (d) $\alpha^{-1}(a)$ is already defined as something other than $i$; (e) $\beta^{-1}(b)$ is already defined as something other than $j$; or (f) $\gamma^{-1}(L[a,b])$ is already defined as something other than $L[i,j]$. In addition, $i$ and $\alpha(i)$ and $j$ and $\beta(j)$ respectively require the same row and column invariants implied by the invariants (as in Sections~\ref{SEI} and \ref{square}).

\subsection{Graph theoretic methods}

A \emph{vertex invariant} is a property of the vertices of graphs which is preserved under automorphisms. In this way, a vertex invariant partitions the vertices such that no automorphism of the graph permutes vertices from one part to another.  In addition to being able to list the automorphisms of the graph, {\sc Nauty} \cite{McKay201494} (available from \url{http://pallini.di.uniroma1.it/})  also outputs a function called \texttt{orbits}, which gives the finest possible vertex invariant partition.

In the graphs we study, the vertex set is or contains $\mathrm{Ent}(L)$ for some partial Latin rectangle $L\in\mathrm{PLR}(r,s,n)$.  Entry invariants can be used to color these vertices.  In general, this coloring does not correspond to a vertex invariant, as it might not be preserved under all automorphisms of the graph, but the automorphisms which have been eliminated by coloring the vertices do not correspond to autotopisms of $L$.

\subsubsection{Adapted McKay, Meynert, and Myrvold method}

We define the vertex-colored graph $G_L$ by
\begin{align*}
V(G_L) := \mathrm{Ent}(L) {}\cup{} & \{R_i \colon\, i \in [r] \text{ and row $i$ of $L$ is non-empty}\} \\
 {}\cup{} & \{S_j \colon\, j \in [s] \text{ and column $j$ of $L$ is non-empty}\} \\
 {}\cup{} & \{N_k \colon\, k \in [n] \text{ and symbol $k$ occurs in $L$}\}
\end{align*}
where each of the four subsets, $\mathrm{Ent}(L)$, $\{R_i\}$, $\{S_j\}$, and $\{N_k\}$, are assigned a distinct vertex color, and edge set
$$
E(G_L) := \{e R_i, e S_j, e N_{L[i,j]} \colon\, e:=(i,j,L[i,j]) \in \mathrm{Ent}(L)\}.
$$
By restricting $\mathrm{Aut}(G_L)$ to $\mathrm{Ent}(L)$, we determine $\mathrm{Atop}(L)$ up to permutations of empty rows and columns, and unused symbols.

\subsubsection{Bipartite graph method}

Pairs of permutations $(\alpha,\beta) \in S_r \times S_s$ act on the set of $r \times s$ $(0,1)$-matrices by permuting the rows by $\alpha$ and the columns by $\beta$; we define an \emph{autotopism} as a stabilizer under this group action.  Corresponding to $L \in \mathrm{PLR}(r,s,n)$, we define a $(0,1)$-matrix $M=M[i,j]_{r \times s}$ with $M[i,j]=1$ if and only if $L[i,j]$ is defined.  The group $\mathrm{Atop}(L)$ is thus isomorphic to a subgroup of the autotopism group of $M$.  Moreover, the matrix $M$ can be interpreted as the biadjacency matrix of a bipartite graph $B_M$.  This gives rise to two ways of computing $\mathrm{Aut}(L)$:
\begin{enumerate}
  \item Use {\sc Nauty} to compute the automorphism group of $B_M$.  This determines the autotopism group of $M$, from which we check each $(\alpha,\beta)$ in the autotopism group of $M$ to see if it defines an autotopism of $L$ through Lemma~\ref{lm:atop_cond}.
  \item Use {\sc Nauty} to compute the \texttt{orbits} of $B_M$.  This determines the finest possible partitions of $[r]$ and $[s]$ for use as row and column invariants, respectively.  We then apply one of the backtracking methods of Section~\ref{se:backtrack} with the added condition that \texttt{orbits} are preserved.
\end{enumerate}

\subsubsection{Partial Latin rectangle graph method}

We define the edge-colored graph $\Gamma_L$ with vertex set $V(\Gamma_L) := \mathrm{Ent}(L)$ and:
\begin{itemize}
 \item green edges between two distinct entries if they share a row,
 \item orange edges between two distinct entries if they share a column, and
 \item purple edges between two distinct entries if they share a symbol.
\end{itemize}
It is known \cite{FalSto16} that any automorphism of $\Gamma_L$ is equivalent to an autotopism of $L$. The vertex set of $\Gamma_L$ is $\mathrm{Ent}(L)$, which implies the automorphisms of $\Gamma_L$ act on the entry set of $L$, so they are determined directly from {\sc Nauty}'s output.

However, {\sc Nauty} does not allow its input graphs to have edge colors; we consider two ways to overcome this.  The first one is to simply ignore the edge colors, then check each automorphism of $\Gamma_L$ returned by {\sc Nauty} to see if it corresponds to an autotopism of $L$.  The second method is to define a new non-edge-colored graph $\overline{\Gamma_L}$, whose automorphisms correspond to the automorphisms of the edge-colored $\Gamma_L$.  It has vertex set

%RAUL + DANI: This original notation seems not to be clear.
%\[
%V(\overline{\Gamma_L}) := \mathrm{Ent}(L) \cup \{e^{(\mathrm{row})}: e \in \mathrm{Ent}(L)\} \cup \{e^{(\mathrm{col})}: e \in \mathrm{Ent}(L)\} \cup \{e^{(\mathrm{sym})}: e \in \mathrm{Ent}(L)\},\]
%where the unions above indicate distinct vertex colors, and edge set
%\begin{align*}
%E(\overline{\Gamma_L}) := \phantom{{}\cup{}} & \{ee^{(\mathrm{row})}, ee^{(\mathrm{col})}, ee^{(\mathrm{sym})}\colon\, e \in \mathrm{Ent}(L)\} \\
%{}\cup{} & \{e_1^{(\mathrm{row})}e_2^{(\mathrm{row})}\colon\, \text{distinct entries } e_1 \text{ and } e_2 \text{ belong to the same row}\} \\
%{}\cup{} & \{e_1^{(\mathrm{col})}e_2^{(\mathrm{col})}\colon\, \text{distinct entries } e_1 \text{ and } e_2 \text{ belong to the same column}\} \\
%{}\cup{} & \{e_1^{(\mathrm{sym})}e_2^{(\mathrm{sym})}\colon\, \text{distinct entries } e_1 \text{ and } e_2 \text{ share the same symbol}\}.
%\end{align*}

\[
V(\overline{\Gamma_L}) := \mathrm{Ent}(L) \cup \{R_e: e \in \mathrm{Ent}(L)\} \cup \{S_e: e \in \mathrm{Ent}(L)\} \cup \{N_e: e \in \mathrm{Ent}(L)\},\]
where the unions above indicate distinct vertex colors, and edge set
\begin{align*}
E(\overline{\Gamma_L}) := \phantom{{}\cup{}} & \{eR_e, eS_e, eN_e\colon\, e \in \mathrm{Ent}(L)\} \\
{}\cup{} & \{R_{e_1}R_{e_2}\colon\, \text{distinct entries } e_1 \text{ and } e_2 \text{ belong to the same row}\} \\
{}\cup{} & \{S_{e_1}S_{e_2}\colon\, \text{distinct entries } e_1 \text{ and } e_2 \text{ belong to the same column}\} \\
{}\cup{} & \{N_{e_1}N_{e_2}\colon\, \text{distinct entries } e_1 \text{ and } e_2 \text{ share the same symbol}\}.
\end{align*}
The graph $\overline{\Gamma_L}$ has more vertices than $\Gamma_L$, so {\sc Nauty} takes longer, but the automorphism group of $\overline{\Gamma_L}$, when restricted to $\mathrm{Ent}(L)$ gives the autotopism group of $L$ directly. Entry invariants can be used to color the vertices in $\mathrm{Ent}(L)$ in both $\Gamma_L$ and $\overline{\Gamma_L}$.

\subsubsection{Rook's graph method}

We define the edge-colored graph $\Xi_L$ from the partial Latin rectangle graph $\Gamma_L$ by deleting the purple edges, thereby forming an induced subgraph of the rook's graph $K_r\times K_s$.  Any autotopism of $L$ is equivalent to an automorphism of $\Xi_L$.  The converse does not necessarily hold, so after {\sc Nauty} has been called, we check each automorphism of $\Xi_L$ to see if it corresponds to an autotopism of $L$.  Differing from partial Latin rectangle graphs, we define $\overline{\Xi_L}$ with
\[
V(\overline{\Xi_L}) := \mathrm{Ent}(L) \cup \{R_e: e \in \mathrm{Ent}(L)\} \cup \{S_e: e \in \mathrm{Ent}(L)\} \cup \{k \in [n]\colon\, \text{symbol } k \text{ appears in } L\},
\]
where the unions above indicate distinct vertex colors, and edge set
\begin{align*}
E(\overline{\Xi_L}) := \phantom{{}\cup{}} & \{eR_e, eS_e\colon\, e \in \mathrm{Ent}(L)\} \\
{}\cup{} & \{R_{e_1}R_{e_2}\colon\, \text{distinct entries } e_1 \text{ and } e_2 \text{ belong to the same row}\} \\
{}\cup{} & \{S_{e_1}S_{e_2}\colon\, \text{distinct entries } e_1 \text{ and } e_2 \text{ belong to the same column}\} \\
{}\cup{} & \{R_e k\colon\, e \in \mathrm{Ent}(L) \text{ and } e \text{ has symbol } k \} \\
{}\cup{} & \{S_e k\colon\, e \in \mathrm{Ent}(L) \text{ and } e \text{ has symbol } k \}.
\end{align*}
It is also possible to use entry invariants to color the vertices in $E(L)$ in both $\Xi_L$ and $\overline{\Xi_L}$.

\section{Experimental results}

We compare run times of all the  methods for computing autotopism groups introduced in Section 3. The experimental platform has an Intel Core i7-4700MQ, 2.40GHz processor (4 physical cores) with 7895MiB RAM running Ubuntu 16.04 LTS (64-bit). In each experiment, we calculate the average run time over $10000$ random $\mathrm{PLR}(r,s,n)$s on the following two randomly generated sets, for some fixed parameters $r$, $s$, and $n$:
\begin{itemize}
    \item In \emph{PLR set A}, we begin with an empty $\mathrm{PLR}(r,s,n)$ and attempt $x$ times to add an entry chosen uniformly at random from $[r] \times [s] \times [n]$.  If a clash arises (that is, if the cell is full, or adding the entry would introduce a repeated symbol in a row or column), we do nothing.  A random partial Latin rectangle generated in this way could have anywhere from $1$ to $rs$ entries, with the expected number of entries growing as the number of attempts $x$ increases.
    \item In \emph{PLR set B}, we generate a random Latin square of order $n\geq \max\{r,s\}$ by using the Jacobson and Matthews method \cite{JacobsonMatthews1996}, then delete the last $n-r$ rows and last $n-s$ columns to obtain an $(rs)$-entry $\mathrm{PLR}(r,s,n)$, that is, a Latin rectangle.  We then delete random entries until we obtain a $\mathrm{PLR}(r,s,n)$ with $x$ entries.
\end{itemize}

Figure~\ref{fi:pilot} in the Appendix plots the experimental run times for random partial Latin rectangles in $\mathrm{PLR}(5,5,5)$. We observe how backtracking methods perform significantly better when invariants are used. Thus, for instance, the run time for entrywise backtracking without the use of entry invariants seems to increase exponentially as the number of entries increases. This fact is not so clear for graph theoretic methods, which show unpredictable performance depending on the input.

\begin{figure}[ht]
\centering
\begin{tabular}{c}
\includegraphics[width=17cm]{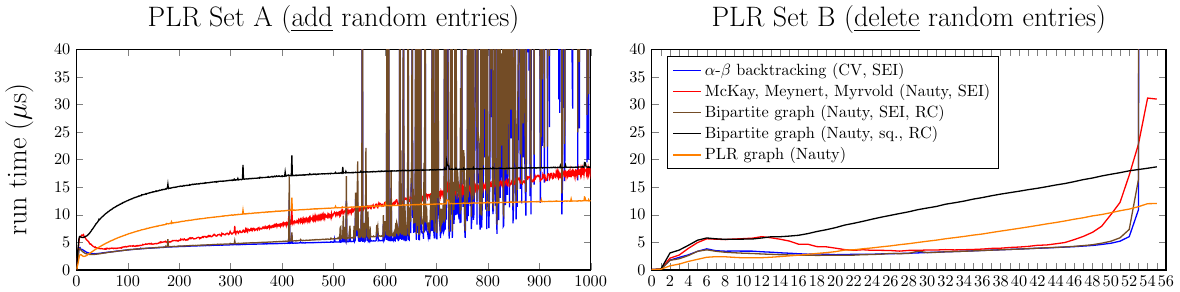}\\
$(r,s,n)=(7,8,9)$\\
\includegraphics[width=17cm]{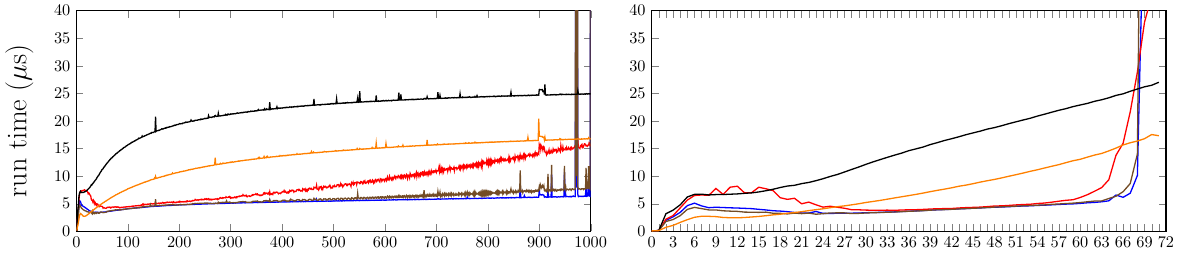}\\
$(r,s,n)=(8,9,10)$\\
\includegraphics[width=17cm]{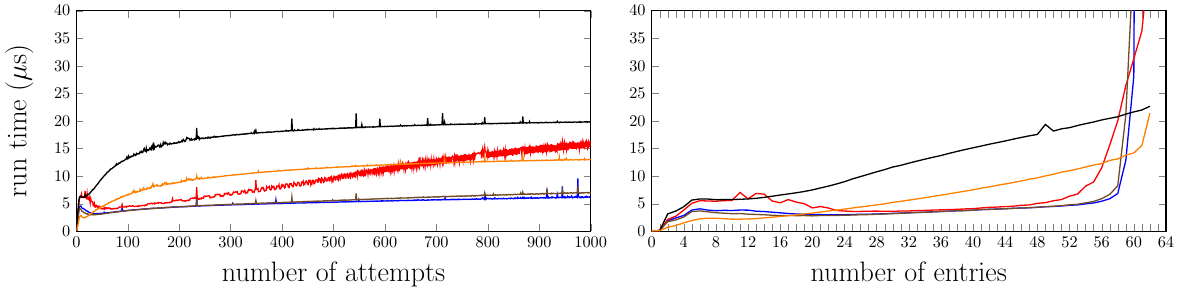}\\
$(r,s,n)=(8,8,8)$\\
\end{tabular}
\caption{The average run time of the remaining methods. For PLR set B, the $63$- and $64$-entry data points for $(r,s,n)=(8,8,8)$ and the $72$-entry data points for $(r,s,n)=(8,9,10)$ are omitted as they took too long to compute.}\label{fi:remain}
\end{figure}

Results for higher orders are shown in Figure~\ref{fi:remain}, where, for PLR Set A, we see fluctuating run times for the alpha-beta backtracking method and bipartite graph method using strong entry invariants. It indicates a high sensitivity to the number of entries.  Therefore, if we compute autotopism groups for partial Latin rectangles where entries are added randomly \cite{Stones13}, we should be aware of the problem that partial Latin rectangles with many entries may drastically reduce the software's performance.

\begin{figure}[htp]
\centering
\begin{tabular}{c}
\includegraphics[width=17cm]{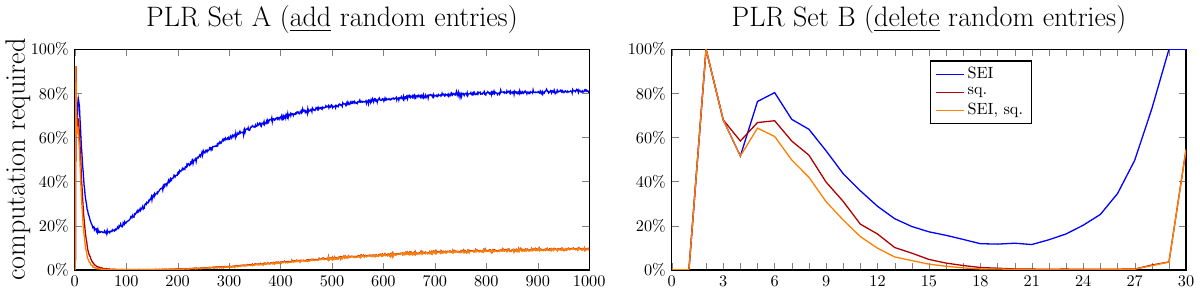}\\
$(r,s,n)=(5,6,7)$\\
\includegraphics[width=17cm]{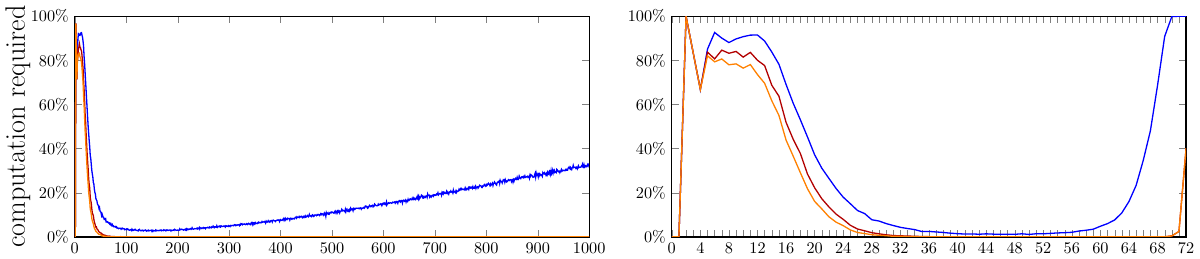}\\
$(r,s,n)=(8,9,10)$\\
\includegraphics[width=17cm]{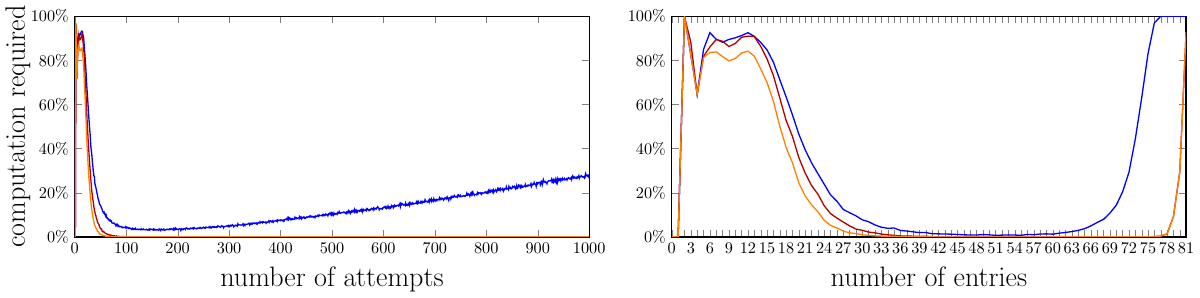}\\
$(r,s,n)=(5,6,7)$\\
\end{tabular}
\caption{The proportion of time (over $10000$ samples) computation is required to compute the autotopism group despite using strong entry invariants (SEI) and/or square invariants (sq.).}\label{fi:comp_required}
\end{figure}

When using entry invariants, sometimes we can deduce that the autotopism group is trivial (up to the autotopisms counted in Lemma~\ref{lm:atop_card}).  In these cases, further computation is not required.  Figure~\ref{fi:comp_required} plots the proportion of the time that further computation is required for three example parameter vectors $(r,s,n)$. We find that entry invariants are ineffective when a partial Latin rectangle has very few entries (likely because it has non-trivial autotopisms).  Further, while strong entry invariants are more efficient to compute than square invariants, Figure~\ref{fi:comp_required} shows that square invariants are more useful at reducing subsequent computation when there are many entries.

\begin{figure}[htp]
\centering

\includegraphics[width=0.6\textwidth]{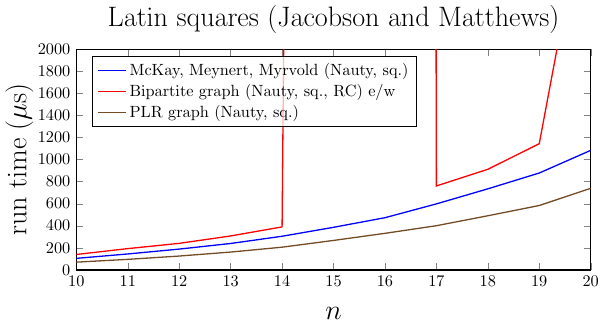}

\caption{The average run time of the MMM method, the bipartite graph method using entrywise backtracking, and the PLR graph method, each using square entry invariants, for random Latin squares of varying order $n$ (here, $r=s=n$).}\label{fi:square}
\end{figure}

\begin{figure}[htp]
\centering

\includegraphics[width=0.6\textwidth]{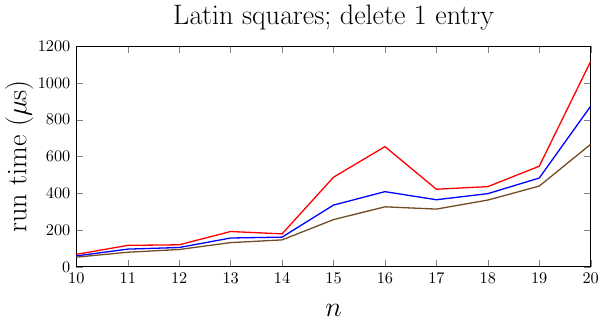}

\includegraphics[width=0.6\textwidth]{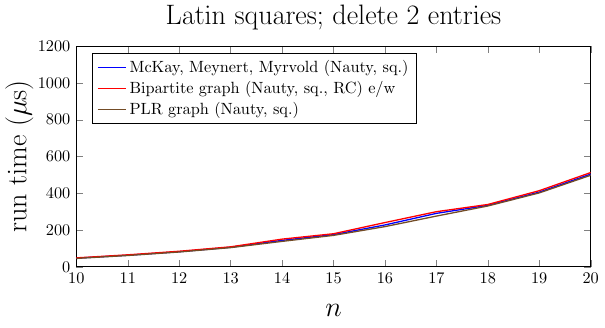}

\caption{The average run time of the MMM method, the bipartite graph method using entrywise backtracking, and the PLR graph method, each using square entry invariants, for random Latin squares of varying order $n$ after deleting either one entry (top plot) or two entries (bottom plot).}\label{fi:square_del}
\end{figure}

Table~\ref{ta:MMMPLR} lists the average run times for the McKay, Meynert, and Myrvold (MMM) method, the bipartite graph method, and the partial Latin rectangle (PLR) graph method when computing autotopism groups for partial Latin rectangles with at most two empty cells.  We see that the PLR graph method outperforms the other methods in every case.  Also clear from Table~\ref{ta:MMMPLR} is the extreme difference between the case $r=s=n$ (that is, Latin squares) and the other cases.  Without entry invariants, the $(17,18,19)$ data set was processed much faster than the $(7,7,7)$; these matrices have 306 cells vs.\ 49 cells, respectively.  In the case of Latin squares, the use of a well-chosen entry invariant is essential to both the MMM and PLR graph methods.  Further experimental run times for Latin squares are given in Figure~\ref{fi:square}, where we see that the PLR graph method consistently outperforms the other methods, and that the bipartite graph method performance is unstable.  Figure~\ref{fi:square_del} includes run times for Latin squares after we delete one or two entries, where the instability is no longer present.

Table~\ref{ta:MMMPLR} shows significantly worse performance when $(r,s,n)=(17,18,19)$ using the MMM or the PLR graph methods with strong entry invariants.  In this setting there is a ``trap'': while there are many entries, there are few distinct strong entry invariants, so coloring the vertices for input into {\sc Nauty} involves sorting a long repetitive list.  % It turns out it is more efficient to leave the vertices uncolored.

\begin{table}[htp]
\centering
\begin{tabular}{|r|rrr|rrr|}
\hline
& \multicolumn{6}{c|}{run time ($\mu$s)} \\
$(r,s,n)$ & \multicolumn{3}{c}{$(17,18,19)$} & \multicolumn{3}{c|}{$(7,7,7)$} \\
\cline{2-7}
no.\ entries & \multicolumn{1}{c}{$rs-0$} & \multicolumn{1}{c}{$rs-1$} & \multicolumn{1}{c|}{$rs-2$} & \multicolumn{1}{c}{$rs-0$} & \multicolumn{1}{c}{$rs-1$} & \multicolumn{1}{c|}{$rs-2$} \\
\hline
MMM ({\sc Nauty}) & 446.0 & 203.1 & 190.3 & 11087.8 & 266.4 & 34.1 \\
MMM ({\sc Nauty}, SEI) & 1714.0 & 262.1 & 293.5 & 11092.8 & 227.5 & 38.7 \\
MMM ({\sc Nauty}, sq.) & 586.2 & 301.3 & 294.9 & 55.2 & 32.1 & 26.5 \\
MMM ({\sc Nauty}, SEI, sq.) & 577.6 & 302.8 & 296.9 & 56.6 & 34.9 & 27.4 \\
\hline
Bipartite graph ({\sc Nauty}, sq., RC) & 777.6 & 295.0 & 292.8 & 5312.0 & 160.0 & 23.0 \\
Bipartite graph ({\sc Nauty}, sq., RC) e/w & 839.8 & 305.7 & 298.4 & 128.0 & 38.2 & 29.0 \\
Bipartite graph ({\sc Nauty}, sq., RC) $\alpha$-$\beta$ & 657.7 & 300.3 & 296.6 & 2230.7 & 87.8 & 25.6 \\
\hline
PLR graph ({\sc Nauty}) & \bf 243.4 & \bf 170.5 & \bf 171.3 & 5476.1 & 154.3 & \bf 17.0 \\
PLR graph ({\sc Nauty}, SEI) & 885.4 & 192.3 & 183.7 & 5482.3 & 37.6 & 20.1 \\
PLR graph ({\sc Nauty}, sq.) & 438.1 & 299.6 & 296.4 & \bf 52.0 & \bf 28.8 & 24.4 \\
PLR graph ({\sc Nauty}, SEI, sq.) & 440.7 & 302.2 & 299.0 & 52.8 & 29.2 & 25.9 \\
\hline
\end{tabular}
\caption{Average run times of the MMM, bipartite graph, and PLR graph methods for partial Latin rectangles with few empty cells (PLR set B).  The best run times are in bold.}\label{ta:MMMPLR}
\end{table}

\section{Conclusions and further work}

From this work, we identify the design goals we should have in mind for software for computing autotopism groups of partial Latin rectangles (PLRs):
\begin{enumerate}
 \item Several methods showed unpredictable performance depending on its input; we should choose a method which has stable performance.
 \item We should mathematically account for PLRs with very few entries, which can have many autotopisms.
 \item A cheap-to-compute entry invariant (such as the strong entry invariant) is useful for PLRs with an intermediate number of entries.
 \item We should use sophisticated entry invariants for PLRs with a large number of entries.
 \item It might be possible to identify PLRs with large (or even transitive) autotopism groups, where we could e.g.\ offload these problematic cases onto e.g.\ a graphics processing unit (GPU) for separate computation.
\end{enumerate}

The PLR graph method offers reasonable and stable performance under most conditions.  It could be improved by
\begin{enumerate}
    \item[a)] finding mathematical conditions on when automorphisms of non-edge-colored partial Latin rectangle graphs are autotopisms of the corresponding partial Latin rectangle, 
    \item[b)] developing software for edge-colored graph automorphism ({\sc Bliss} \cite{JunttilaKaski2007}, available from \url{http://www.tcs.hut.fi/Software/bliss/}, also does not incorporate edge colors), or 
    \item[c)] developing purpose-built software which adapts the individualization/refinement approach of {\sc Nauty} and {\sc Bliss} to the problem of computing autotopisms of partial Latin rectangles.
\end{enumerate}

Further work on the development of efficient algorithms for the computational processes that have been proposed throughout the paper, along with a study of the reasons that give rise to better or worse experimental results, is also required.

% becky: they do, and that's part of my new paper---this is no longer "future research"
% Recently $3$-plexes of Latin squares (a special case of $k$-plexes \cite{Wanless2002}, which is a special type of partial Latin square) have been used to develop erasure codes \cite{YiEtAl}.  This suggests that studying $k$-plexes with autotopisms may have applications in the study of erasure codes.

\subsection*{Acknowledgements}
Marbach and Stones's work is partially supported by NSF of China (61602266, 61872201), the Science and Technology Development Plan of Tianjin (17JCYBJC15300, 16JCYBJC41900), and the Fundamental Research Funds for the Central Universities and SAFEA: Overseas Young Talents in Cultural and Educational Sector.  Stones is also supported by the Thousand Youth Talents Plan in Tianjin.

Stones was supported by her NSFC Research Fellowship for International Young Scientists (grant numbers: 11450110409, 11550110491), and the Thousand Youth Talents Plan in Tianjin.

Falc\'on's work is partially supported by the research project FQM-016 from Junta de Andaluc\'\i a, and the Departmental Research Budget of the Department of Applied Mathematics I of the University of Seville.

The authors thank Zhuanhao Wu for assistance coding.

\bibliographystyle{wileyNJD-AMA}
\bibliography{atop}

\appendix
\section{Experimental plots}

\begin{figure}[h!]
\centering
\includegraphics[width=17cm]{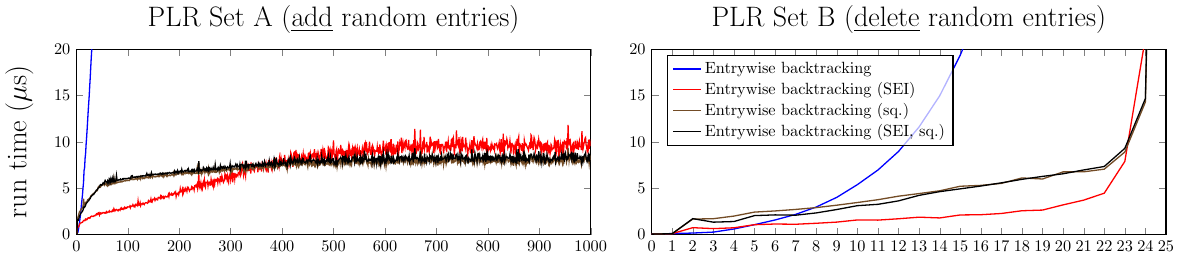}
\includegraphics[width=17cm]{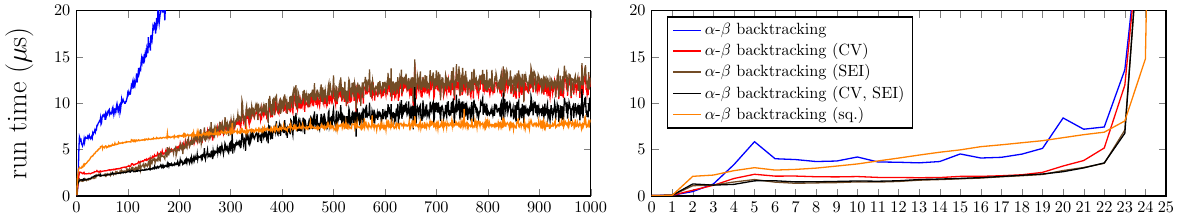}
\includegraphics[width=17cm]{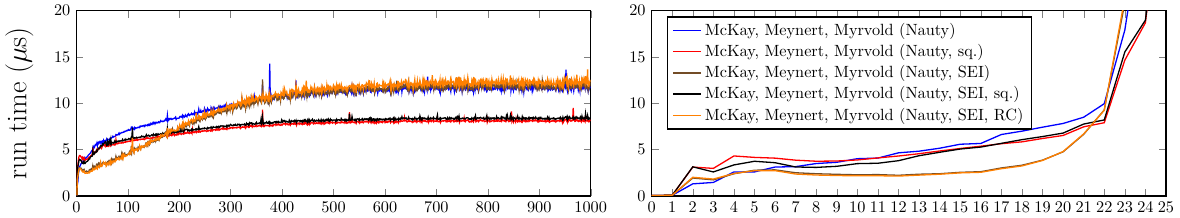}
\includegraphics[width=17cm]{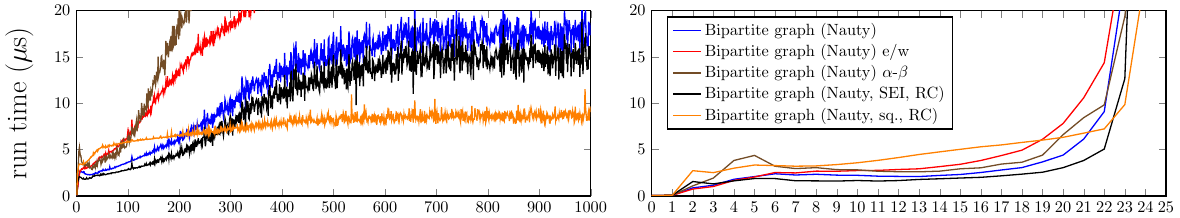}
\includegraphics[width=17cm]{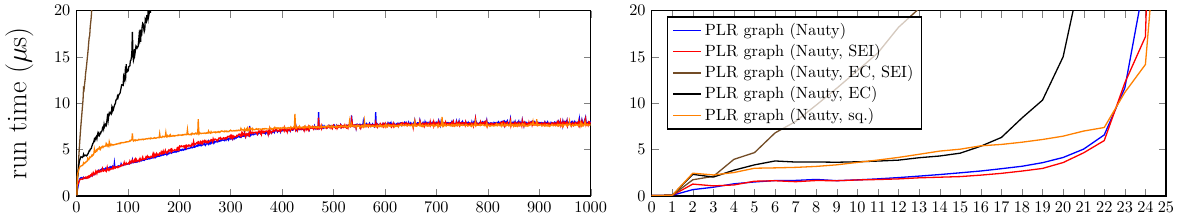}
\includegraphics[width=17cm]{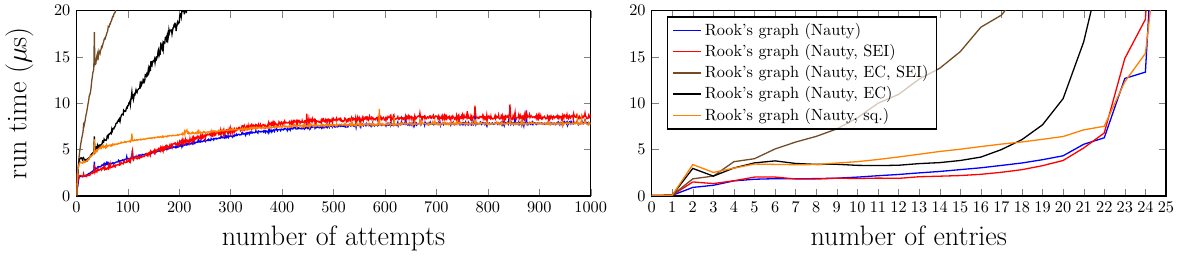}

\caption{The average run time of various methods used to compute the autotopism group, with or without the use of strong entry invariants (SEI) and/or square invariants (sq.); parameters $(r,s,n)=(5,5,5)$.  Other acronyms: RC = ``row and column invariants''; e/w = ``entrywise''; EC = ``edge colors''; CV = ``column vectors''.  Some experimental results are omitted for space reasons.}\label{fi:pilot}
\end{figure}

\end{document}